\documentclass[12pt]{article}
\usepackage{amsmath,amsthm,amssymb}
\usepackage{mathrsfs}
\usepackage[titletoc]{appendix}
\usepackage{bbm}
\usepackage{float}
\usepackage{tikz}
\begin{document}
\input amssym.def
\newcommand{\singlespace}{
    \renewcommand{\baselinestretch}{1}
\large\normalsize}
\newcommand{\doublespace}{
   \renewcommand{\baselinestretch}{1.2}
   \large\normalsize}
\renewcommand{\theequation}{\thesection.\arabic{equation}}

\setcounter{equation}{0}
\def \ten#1{_{{}_{\scriptstyle#1}}}
\def \Z{\Bbb Z}
\def \C{\Bbb C}
\def \R{\Bbb R}
\def \Q{\Bbb Q}
\def \N{\Bbb N}
\def \l{\lambda}
\def \V{V^{\natural}}
\def \wt{{\rm wt}}
\def \tr{{\rm tr}}
\def \Res{{\rm Res}}
\def \End{{\rm End}}
\def \Aut{{\rm Aut}}
\def \mod{{\rm mod}}
\def \Hom{{\rm Hom}}
\def \im{{\rm im}}
\def \<{\langle}
\def \>{\rangle}
\def \w{\omega}
\def \c{{\tilde{c}}}
\def \o{\omega}
\def \t{\tau }
\def \ch{{\rm ch}}
\def \a{\alpha }
\def \b{\beta}
\def \e{\epsilon }
\def \la{\lambda }
\def \om{\omega }
\def \O{\Omega}
\def \qed{\mbox{ $\square$}}
\def \pf{\noindent {\bf Proof: \,}}
\def \voa{vertex operator algebra\ }
\def \voas{vertex operator algebras\ }
\def \p{\partial}
\def \1{{\bf 1}}
\def \ll{{\tilde{\lambda}}}
\def \H{{\bf H}}
\def \F{{\bf F}}
\def \h{{\frak h}}
\def \g{{\frak g}}
\def \rank{{\rm rank}}
\def \({{\rm (}}
\def \){{\rm )}}
\def \Y {\mathcal{Y}}
\def \I {\mathcal{I}}
\def \A {\mathcal{A}}
\def \B {\mathcal {B}}
\def \Cc {\mathcal {C}}
\def \H {\mathcal{H}}
\def \M {\mathcal{M}}
\def \V {\mathcal{V}}
\def \O{{\bf O}}
\def \AA{{\bf A}}
\singlespace
\newtheorem{thm}{Theorem}[section]
\newtheorem{prop}[thm]{Proposition}
\newtheorem{lem}[thm]{Lemma}
\newtheorem{cor}[thm]{Corollary}
\newtheorem{rem}[thm]{Remark}
\newtheorem*{CPM}{Theorem}
\newtheorem{definition}[thm]{Definition}

\begin{center}
{\Large {\bf  Representations of vertex operator algebras and bimodules}} \\

\vspace{0.5cm} Chongying Dong\footnote{Supported by a NSF grant.}
\\
Department of Mathematics\\ University of
California\\ Santa Cruz, CA 95064 \\
Li Ren\footnote{Supported in part by China Postdoctor grant 2012M521688.}\\
 School of Mathematics \\ Sichuan University\\
Chengdu 610064, China\\
\end{center}
\hspace{1.5 cm}

\begin{abstract}
For a vertex operator algebra $V,$  a $V$-module $M$ and a nonnegative integer $n,$ an $A_n(V)$-bimodule
$\AA_n(M)$ is constructed and studied.  The connection between $\AA_n(M)$ and intertwining operators
are discussed. Moreover, the $A_n(V)$-bimodule $A_{t,s}(V)$ is a quotient of $\AA_n(V)$ for all $s,t\leq n.$
In the case that $V$ is rational, $\AA_n(M)$ for irreducible $V$-module $M$ is given explicitly.

\end{abstract}

\section{Introduction}

The associative algebra $A_n(V)$ for any vertex operator algebra $V$ and a nonnegative integer $n$ was constructed in
\cite{DLM3} such that $A_{n-1}(V)$ is a quotient of $A_n(V)$ and  there is a bijection between irreducible
admissible $V$-modules and irreducible $A_n(V)$-modules which cannot factor through $A_{n-1}(V).$
Moreover, $V$ is rational if and only if $A_n(V)$ are finite dimensional semisimple associative algebras for all $n.$
In the case $n=0,$ $A_0(V)$ is exactly the algebra $A(V)$ investigated in \cite{Z}. An $A(V)$-bimodule $A(M)$
for any $V$-module $M$ was also introduced in \cite{FZ} to deal with the intertwining operators and fusion rules.

Following the $A(M)$-theory from \cite{FZ} we construct a sequence of $A_n(V)$-bimodules $\AA_n(M)$
for any $V$-module $M$ and a nonnegative integer $n$ such that $\AA_{n-1}(M)$ is a quotient of $\AA_n(M)$ and
$\AA_0(M)=A(M).$ Moreover, the $A_n(V)$-bimodule $A_{t,s}(V)$ for $s,t\leq n$  defined in \cite{DJ} is a quotient of $\AA_n(M).$ It is established that if $V$ is rational then there is a linear isomorphism from
the space  ${\cal I}_{M^i\,M^j}^{M^k}$ of intertwining operators to
$$\Hom_{A_n(V)}(\AA_n(M^i)\otimes_{A_n(V)}M^j(s), M^k(t))$$
for $s,t\leq n$ where $M^q=\oplus_{m\geq 0}M^q(m)$ $(q=i,j,k)$ are the irreducible $V$-module such that $M^q(0)\ne 0.$ This result is a generalization of that obtained in \cite{FZ} when $n=0.$ 
The bimodule structure of $\AA_n(M)$ for any irreducible module for rational $V$ is given explicitly. 

One can regard $\AA_n$ as a functor from the $V$-module category to the $A_n(V)$-bimodule category. One important property of the functor $\AA_0=A$ is that $A$ respects to the tensor product at least for rational vertex operator algebra. That is, $A(M\boxtimes N)=A(M)\otimes_{A(V)}A(N)$ for any $V$-modules
$M,N$ where $M\boxtimes N$ is the tensor product of $V$-modules as studied in \cite{HL1, HL2, H}. 
Unfortunately, this is not true for general $n.$ This can be seen clearly from the explicit $A_n(V)$-bimodule structure of $\AA_n(M).$ 

Another interesting result about $\AA_n(M)$ is the relation between $\AA_n(M)$ and $\AA_n(M')$ where $M'$ is the contragredient module of $M$ as defined in \cite{FHL}. It is well known from the bimodule theory that
$\AA_n(M)^*$ is also an $A_n(V)$-bimodule in an obvious way. We show that if $V$ is rational and $C_2$-cofinite then
$\AA_n(M)^*$ is isomorphic to $\AA_n(M')$ for any irreducible $V$-module $M.$ The proof involves a relation
on the fusion matrices associated to $M$ and $M'$ \cite{DJX}. This explains why we can only prove the isomorphism
between $\AA_n(M)^*$ and $\AA_n(M')$ under rationality and $C_2$-cofiniteness assumptions. We certainly believe that this result is true in general as long as $\AA_n(M)$ is finite dimensional. A proof of this result without using the fusion matrices will be important and useful.

Note that our $\AA_n(M)$ is different from $A_n(M)$  defined in \cite{HY} where $A_n(M)=M/O_n(M)$ and
$O_n(M)$ also contains $(L(-1)+L(0))M.$ From the connection between intertwining operators and $\AA_n(M)$ discussed below it seems that $(L(-1)+L(0))M$ should not be a subspace of $\O_n(M)$ in our consideration. 

One of the important motivations for constructing $A_n(V)$-bimodule $\AA_n(M)$ is to study the extension of
rational vertex operator algebras. It is a well known conjecture that if  $V$ is a rational vertex operator
algebra then any extension $U$ of $V$ is also rational. It is expected that the $\AA_n(M)$-theory will play roles in proving this conjecture.

There are associative algebras $A_{g,n}(V)$ associated to an automorphism $g$ of $V$ of finite order and
$n\in \frac{1}{o(g)}\Z_+$ \cite{DLM4}, \cite{MT}.  One could  construct $A_{g,n}(V)$-bimodule
 $A_{g,n}(M)$ for a $g$-twisted $V$-module $M$ following the ideals of this paper. 

The paper is organized as follows. We review the construction of associative algebras $A_n(V)$ and relevant results
from  \cite{DLM3} in Section 2. The construction of $\AA_n(M)$ is given in Section 3. We also show how the identity map on $M$ induces an $A_n(V)$-bimodule epimorphism from $\AA_n(M)$ to $\AA_{n-1}(M).$ Section 4 is devoted to the study of relation between $\AA_n(M)$ and intertwining operators. As in \cite{FZ} and \cite{L2} we argue how the map from ${\cal I}_{M^i\,M^j}^{M^k}$ to $\Hom_{A_n(V)}(\AA_n(M^i)\otimes_{A_n(V)}M^j(s), M^k(t))$ by sending
$I\in {\cal I}_{M^i\,M^j}^{M^k}$ to $I_{t,s}$ which maps $M^j(s)$ to $M^k(t)$ induces a bijection from ${\cal I}_{M^i\,M^j}^{M^k}$ to
$$\Hom_{A_n(V)}(\AA_n(M^i)\otimes_{A_n(V)}M^j(s), M^k(t))$$
if $V$ is rational. Various properties of $\AA_n(M)$ are discussed. In Section 5 we investigate the relation between
$\AA_n(M)^*$ and $\AA_n(M').$

\section{$A_n(V)$ construction}
\setcounter{equation}{0}

This section is a review of the associative algebra $A_n(V)$ and  related results from \cite{DLM3}. Also see \cite{Z} .

Let $V=(V,Y,{\bold 1},\omega)$ be a vertex operator algebra \cite{B}, \cite{FLM}, \cite{LL}.
 We first recall different notions of modules from \cite{DLM1, FLM, Z}. A {\em weak V-module} $M$ is a vector space
equipped with a linear map
\begin{equation*}
\begin{split}
Y_M(\cdot, z) : &V \to (\End M)[[z, z^{-1}]]\\
&v  \mapsto Y_M(v, z) =\sum_{n\in\Z}v_nz^{-n-1}~~~~~~(v_n \in \hbox{End}M)
\end{split}
\end{equation*}
which satisfies the following conditions for $u \in V,$ $v \in V,$ $w \in M$ and $ n \in \Z,$
\begin{equation*}
\begin{split}
&u_nw = 0 \hbox{ for } n \gg 0;\\
& Y_M(\1, z) = \hbox{id}_M;\\
z_0^{-1}\delta(&\frac{z_1-z_2}{z_0})Y_M(u, z_1)Y_M(v, z_2)-z_0^{-1}\delta(\frac{z_2-z_1}{-z_0})Y_M(v, z_2)Y_M(u, z_1)\\
&= z_2^{-1}\delta(\frac{z_1-z_0}{z_2})Y_M(Y (u, z_0)v, z_2).
\end{split}
\end{equation*}

An ({\em ordinary}) $V$-module is a weak $V$-module $M$ which carries a $\C$-grading induced by the spectrum
of $L(0)$ where $L(0)$ is a component operator of
$$Y_M(\omega, z) =\sum_{n\in\Z}L(n)z^{-n-2}.$$
That is, $M =\oplus _{\l\in \C}M_{\l}$ where $M_\l = \{w\in M|L(0)w = \l w\}$. Moreover one requires that $M_\l$ is
finite dimensional and for fixed $\l$, $M_{n+\l}= 0$ for all small enough integers $n$.

An {\em admissible} $V$-module is a weak $V$-module $M$ which carries a $\Z_+$-grading
$M =\oplus_{n\in\Z_+} M(n)$ that satisfies the following
$$v_mM(n) \subset M(n+\wt v-m-1)$$
for homogeneous $v\in V.$ It is easy to show that any {\em ordinary} module is {\em admissible}.

For an {\em ordinary} $V$-module $M=\bigoplus_{\lambda\in \C}M_{\lambda},$ the contragredient module $M'$
is defined in \cite{FHL} as follows:
\begin{equation*}
M'=\bigoplus_{n\in \lambda}M_{\lambda}^{*},
\end{equation*}
where $M_{\lambda}^*$ is the dual space of $M_{\lambda}.$ The vertex operator
$Y_{M'}(a,z)$ is defined for $a\in V$ via
\begin{eqnarray*}
\langle Y_{M'}(a,z)f,w\rangle= \langle f,Y_M(e^{zL(1)}(-z^{-2})^{L(0)}a,z^{-1})w\rangle,
\end{eqnarray*}
where $\langle f,w\rangle=f(w)$ for $f\in M', w\in M$ is the natural paring $M'\times M\to \C.$

$V$ is called {\em rational} if every admissible $V$-module is completely reducible.  It is proved in \cite{DLM2}
that if $V$ is rational then there are only finitely many irreducible admissible $V$-modules up to
isomorphism and each irreducible admissible $V$-module is ordinary.  Let
$M^0, \ldots,M^p$ be the irreducible modules up to isomorphism with $M^0 = V$. Then there exist $\lambda_i \in\C$
for $i = 0, \ldots, p$ such that
$$M^i = \oplus_{n=0}^\infty M^i_{\l_i+n}$$
 with $M^i_{\l_i}\neq 0$ where $L(0)|_{M^i_{\l_i+n}}= \l_i + n$  and $ n\in \Z_+.$ $\l_i$ is called the {\em conformal weight} of
 $M^i.$  We denote $M^i(n)=M^i_{\l_i+n}.$ Moreover, $\l_i$ and the central charge $c$ are
rational numbers (see \cite{DLM5}). Let $(M^i)'=M^{i^*}$ for some unique $i^*\in\{0,...,p\}.$

$V$ is called $C_2$-{\em cofinite} if $\hbox{dim} V/C_2(V ) <\infty$ where $ C_2(V ) = \< u_{-2}v|u, v \in V \>$ \cite{Z}.

We now define algebra $A_n(V)$ for nonnegative integer $n.$ Let $O_n(V)$ be the linear span of $u\circ_n v$ and $L(-1)u+L(0)u$
where for homogeneous $u\in V$ and $v\in V,$
\begin{equation*}\label{2.1}
u\circ_n v=\Res_{z}Y(u,z)v\frac{(1+z)^{\wt u+n}}{z^{2n+2}}.
\end{equation*}
Also define second product $*_n$ on $V$ for $u$ and $v$ as
above:
\begin{eqnarray*}\label{2.2}
& & u*_nv=\sum_{m=0}^{n}(-1)^m{m+n\choose n}\Res_zY(u,z)\frac{(1+z)^{\wt\,u+n}}{z^{n+m+1}}v\\
& &\ \ \ \ \ \ \ \ =\sum_{m=0}^n\sum_{i=0}^{\infty}(-1)^m
{m+n\choose n}{\wt u+n\choose i}u_{i-m-n-1}v.\nonumber
\end{eqnarray*}
Extend linearly to obtain a bilinear product on $V.$  Set $A_n(V)=V/O_n(V).$

The following theorem summarizes the main results of \cite{DLM3}.
\begin{thm}\label{tha} Let $V$ be a vertex operator algebra and $n$ a nonnegative
integer. Then

(1) $A_n(V)$ is an associative algebra whose product is induced by
$*_n.$

(2) The identity map on $V$ induces an algebra epimorphism from
$A_n(V)$ to $A_{n-1}(V).$

(3) Let $W$ be a weak module and set
$$\Omega_n(W)=\{w\in W|u_mw=0, u\in V, m> \wt u-1+n\}.$$
Then $\Omega_n(W)$ is an $A_n(V)$-module such that $v+O_n(V)$ acts as $o(v)$ where
$o(v)=v_{\wt v-1}$ for homogeneous $v$ and extend linearly to entire $V.$

(4) Let $M=\oplus_{m=0}^{\infty}M(m)$ be an admissible $V$-module.
Then each $M(m)$ for $m\leq n$ is
an $A_n(V)$-submodule of $\Omega_n(W).$ Furthermore,
$M$ is irreducible if and only if each $M(n)$ is an irreducible
$A_n(V)$-module for all $n.$

(5) $V$ is rational if and only if $A_n(V)$ are finite dimensional
semisimple algebras for all $n\geq 0.$ In this case
$$A_n(V)=\bigoplus_{i=0}^p\bigoplus_{j=0}^n\End M^i(j)
=\bigoplus_{i=0}^p\bigoplus_{j=0}^nM^i(j)\otimes_{\C}M^{i^*}(j)$$
where $M^i$ for $i=0,...,p$ are the irreducible $V$-modules.
Moreover, $O_{s-1}(V)/O_{s}(V)$ and $A_{s-1}(V)$ are two sided ideals of $A_n(V)$ for $s=1,...,n,$
and
$$O_{s-1}(V)/O_{s}(V)=\bigoplus_{i=0}^p\End M^i(s).$$

(6) The map $M\mapsto M(n)$ gives  one to one correspondence between irreducible
admissible $V$-modules and the irreducible $A_n(V)$-modules which are not $A_{n-1}(V)$-modules.
\end{thm}

The most parts of Theorem \ref{tha} are clear. We give a few words on Theorem \ref{tha} (5). By (1)
the identity map on $V$ induces an algebra epimorphism from $A_n(V)$ to $A_{n-1}(V)$ with kernel
$O_{n-1}(V)/O_n(V).$ Since $A_n(V)$ is semisimple we conclude that both $O_{n-1}(V)/O_n(V)$ and $A_{n-1}(V)$ are the ideals of $A_n(V).$ It is now clear that
$$O_{n}(V)/O_{n-1}(V)=\bigoplus_{i=0}^p\End M^i(n).$$
The decomposition of $O_{s}(V)/O_{s-1}(V)$ for arbitrary $s$ follows immediately.

\section{$A_n(V)$-bimodules $\AA_n(M)$}

Let $V$ be a vertex operator algebra and $M$ be an admissible $V$-module. Motivated by the $A(V)$-bimodule $A(M)$ from \cite{Z} we define and study the $A_n(V)$-bimodule $\AA_n(M)$ for any nonnegative integer $n.$ The construction of $\AA_n(M)$ is largely influenced by the construction of $A_n(V)$ \cite{DLM3} and the intertwining operators \cite{FHL}. See \cite{HY} for a different treatment.

Let $\O_n(M)$ be the linear span of $u\circ_n w$
where for homogeneous $u\in V$ and $w\in M,$
\begin{equation*}\label{2.1m}
u\circ_n w=\Res_{z}Y(u,z)w\frac{(1+z)^{\wt u+n}}{z^{2n+2}}.
\end{equation*}
Also define a left bilinear product $*_n$  for $u\in V$ and $w\in M:$
\begin{eqnarray}\label{2.2m}
& & u*_nw=\sum_{m=0}^{n}(-1)^m{m+n\choose n}\Res_zY(u,z)\frac{(1+z)^{\wt\,u+n}}{z^{n+m+1}}w\\
& &\ \ \ \ \ \ \ \ =\sum_{m=0}^n\sum_{i=0}^{\infty}(-1)^m
{m+n\choose n}{\wt u+n\choose i}u_{i-m-n-1}w,\nonumber
\end{eqnarray}
and a right bilinear  product
\begin{eqnarray}\label{2.3m}
& & w*_nu=\sum_{m=0}^{n}(-1)^n{m+n\choose n}\Res_zY(u,z)\frac{(1+z)^{\wt\,u+m-1}}{z^{n+m+1}}w\\
& &\ \ \ \ \ \ \ \ =\sum_{m=0}^n\sum_{i=0}^{\infty}(-1)^m
{m+n\choose n}{\wt u+m-1\choose i}u_{i-m-n-1}w.\nonumber
\end{eqnarray}
Set $\AA_n(M)=M/\O_n(M).$ In the case $n=0$ the $\AA_0(M)$ is exactly the $A(V)$-bimodule $A(M)$ studied in \cite{FZ}

\begin{rem} We have already mentioned that our $\AA_n(M)$ is different from $A_n(M)$  defined in \cite{HY} where $A_n(M)=M/O_n(M)$ and $O_n(M)$ also contains
$(L(-1)+L(0))M.$
\end{rem}

\begin{rem} In the case $n=0$ it follows immediately from the definitions that $\AA_0(V)=A_0(V)=A(V)$ \cite{Z}.
But it is not clear if this is true in general. It will be established later that $\AA_n(V)=A_n(V)$ for all $n$ if
$V$ is rational.
\end{rem}

\begin{lem}\label{l2.2} (1)
 Assume that $u\in V$ is homogeneous,
$w\in M$ and $m\ge k\ge 0.$ Then
$$\Res_{z}Y(u,z)w\frac{(1+z)^{{\wt}u+n+k}}{z^{2n+2+m}}\in \O_n(M).$$

(2) For homogeneous $u\in V$ and $w\in M,$
$u*_nw-w*_nu-\Res_zY(u,z)w(1+z)^{\wt u-1}\in \O_{n}(M).$
\end{lem}

\pf The proof of (i) is similar to that of Lemma 2.1.2 of [Z]. (ii) follows
from the definitions.
\qed
\begin{lem}\label{l2.3} We have the following containments:

(1) $O_n(V)*_nM\subset \O_n(M), M*_nO_n(V)\subset \O_n(M),$

(2) $V*_n\O_n(M)\subset \O_n(M), \O_n(M)*_nV\subset \O_n(M).$
\end{lem}

\pf The argument of containments $(L(-1)+L(0))V*_nM\subset \O_n(M), M*_n (L(-1)+L(0))V\subset \O_n(M)$
is the similar to that of $(L(-1)+L(0))V*_nV\subset O_n(V), V*_n (L(-1)+L(0))V\subset O_n(V)$
presented in the proof of Lemma 2.2 of \cite{DLM3} using Lemma \ref{2.2} (2). It remains to prove that
$(u\circ_nv)*_nw, w*_n(u\circ_n v), u*_n(v\circ_nw), (u\circ_nw)*_nv\in \O_n(M)$ for $u,v\in V$ and $w\in M.$ Since the proofs are similar we only prove $(u\circ_nv)*_nw\in \O_n(M).$

We have
\begin{eqnarray*}
& & (u\circ_nv)*_nw=\sum_{i\geq 0}{\wt u+n\choose i}(u_{-2n-2+i}v)*_nw\\
& &=\sum_{i\geq 0}\sum_{m=0}^n(-1)^{m}{m+n\choose
n}{\wt u+n\choose i}\Res_{z_2}Y(u_{-2n-2+i}v,z_2)w
\frac{(1+z_2)^{\wt u+\wt v+3n+1-i}}{z_2^{1+m+n}}\\
& &=\sum_{m=0}^n(-1)^{m}{m+n\choose n}\Res_{z_2}\Res_{z_1-z_2}Y(Y(u,z_1-z_2)v,z_2)w\\
& &\ \ \ \ \ \cdot
\frac{(1+z_1)^{\wt u+n}(1+z_2)^{\wt v+2n+1}}{(z_1-z_2)^{2n+2}z_2^{1+m+n}}\\
& &=\sum_{m=0}^n(-1)^{m}{m+n\choose
n}\Res_{z_1}\Res_{z_2}Y(u,z_1)Y(v,z_2)w
\frac{(1+z_1)^{\wt u+n}(1+z_2)^{\wt v+2n+1}}{(z_1-z_2)^{2n+2}z_2^{1+m+n}}\\
& &\ \ \ \ \ -\sum_{m=0}^n(-1)^{m}{m+n\choose
n}\Res_{z_2}\Res_{z_1}Y(v,z_2)Y(u,z_1)w\frac{(1+z_1)^{\wt u+n}(1+z_2)^{\wt v+2n+1}}{(z_1-z_2)^{2n+2}z_2^{1+m+n}}\\
& &=\sum_{m=0}^n\sum_{i\geq 0}(-1)^{m}{m+n\choose
n}{-2n-2\choose i}(-1)^i\\
& &\ \ \ \ \ \Res_{z_1}\Res_{z_2}Y(u,z_1)Y(v,z_2)w\frac{(1+z_1)^{\wt u+n}(1+z_2)^{\wt v+2n+1}}{z_1^{2n+2+i}z_2^{1+m+n-i}}\\
& &\ \ \ \ \ -\sum_{m=0}^n\sum_{i\geq 0}(-1)^{m}{m+n\choose
n}{-2n-2\choose i}(-1)^i\\
& &\ \ \ \ \ \Res_{z_2}\Res_{z_1}Y(v,z_2)Y(u,z_1)w
\frac{(1+z_1)^{\wt u+n}(1+z_2)^{\wt v+2n+1}}{z_1^{-i}
z_2^{3n+3+m+i}}.
\end{eqnarray*}
From Lemma \ref{l2.2} we know that both
$$\Res_{z_1}\Res_{z_2}Y(u,z_1)Y(v,z_2)w\frac{(1+z_1)^{\wt u+n}(1+z_2)^{\wt v+2n+1}}{z_1^{2n+2+i}z_2^{1+m+n-i}}$$
and
$$\Res_{z_2}\Res_{z_1}Y(v,z_2)Y(u,z_1)w\frac{(1+z_1)^{\wt u+n}(1+z_2)^{\wt v+2n+1}}{z_1^{-i}
z_2^{3n+3+m+i}}$$
lie in $\O_n(M).$ The proof is complete. \qed

Recall from \cite{DLM3} that the linear map
$$\phi:  v\mapsto e^{L(1)}(-1)^{L(0)}v$$
induces an anti-isomorphism $A_n(V)$ to itself. We can now establish the following:
\begin{thm}\label{t3.5} Let $M$ be an admissible $V$-module and  $n\geq 0.$

(1) The $\AA_n(M)$ is an $A_n(V)$-bimodule such that the left and right actions of $A_n(V)$ on $\AA_n(M)$
induced from (\ref{2.2m}) and (\ref{2.3m}), respectively.

(2) The identity map on $M$ induces an $A_n(V)$-bimodule epimorphism from $\AA_n(M)$ to $\AA_{n-1}(M)$
if $n\geq 1.$

(3) The map
$$\phi:  w\mapsto e^{L(1)}e^{\pi iL(0)}w$$
induces a linear isomorphism from $\AA_n(M)$ to itself such that
$$\phi(u*_nw)=\phi(w)*_n\phi(u), \phi(w*_nu)=\phi(u)*_n\phi(w)$$
for $u\in V$ and $w\in M.$

(4) If $V$ is rational then both $\O_{s-1}(M)/\O_{s}(M)$ and $\AA_{s-1}(M)$ are the $A_n(V)$-bimodules
for $s=1,...,n$ and
$$\AA_n(M)=\AA_0(M)\bigoplus \bigoplus_{s=1}^n\O_{s-1}(M)/\O_{s}(M).$$
\end{thm}

\pf (1) By Lemma \ref{l2.3} it is good enough to prove the following relations in $\AA_n(M)$ for $u,v\in V$ and $w\in M:$
\begin{eqnarray*}
& &(u*_nw)*_nv=u*_n(w*_nv), \1*_nw=w*_n\1=w\\
& & (u*_nv)*_nw=u*_n(v*_nw), w*_n(u*_nv)=(w*_nu)*_nv.
\end{eqnarray*}
Again the proofs are similar for these relations. We give a detail proof  for $(u*_nw)*_nv=u*_n(w*_nv).$

Using the definition we have
\begin{eqnarray*}
& &u*_n(w*_nv)-(u*_nw)*_nv=\sum_{m_1,m_2=0}^n(-1)^{m_1+n}{m_1+n\choose
n}{m_2+n\choose n}\\
& &\cdot \Res_{z_1}\Res_{z_2}Y(u,z_1)Y(v,z_2)w\frac{(1+z_1)^{\wt\,u+n}(1+z_2)^{\wt v+m_2-1}}{z_1^{n+m_1+1}z_2^{n+m_2+1}}\\
& &-\sum_{m_1,m_2=0}^n(-1)^{m_1+n}{m_1+n\choose
n}{m_2+n\choose n}\\
& &\cdot \Res_{z_2}\Res_{z_1}Y(v,z_2)Y(u,z_1)w\frac{(1+z_1)^{\wt\,u+n}(1+z_2)^{\wt v+m_2-1}}{z_1^{n+m_1+1}z_2^{n+m_2+1}}\\
& &=\sum_{m_1,m_2=0}^n(-1)^{m_1+n}{m_1+n\choose
n}{m_2+n\choose n}\\
& &\cdot \Res_{z_2}\Res_{z_1-z_2}Y(Y(u,z_1-z_2)v,z_2)w\frac{(1+z_1)^{\wt\,u+n}(1+z_2)^{\wt v+m_2-1}}{z_1^{n+m_1+1}z_2^{n+m_2+1}}\\
& &=\sum_{m_1,m_2=0}^n\sum_{i,j\geq 0}(-1)^{m_1+n}{m_1+n\choose
n}{m_2+n\choose n}{\wt u+n\choose i}{-n-m_1-1\choose j}\\
& &\cdot \Res_{z_2}\Res_{z_1-z_2}Y(Y(u,z_1-z_2)v,z_2)w\frac{(1+z_2)^{\wt u+\wt v+n+m_2-1-i}(z_1-z_2)^{i+j}}{z_2^{2n+m_1+m_2+2+j}}\\
& &=\sum_{m_1,m_2=0}^n\sum_{i,j\geq 0}(-1)^{m_1+n}{m_1+n\choose
n}{m_2+n\choose n}{\wt u+n\choose i}{-n-m_1-1\choose j}\\
& &\cdot \Res_{z_2}Y(u_{i+j}v,z_2)w\frac{(1+z_2)^{\wt u+\wt v+n+m_2-1-i}}{z_2^{2n+m_1+m_2+2+j}}.
\end{eqnarray*}
Note that $\wt u_{i+j}v=\wt u+\wt v-i-j-1$ and  $\wt u+\wt v+n+m_2-1-i=\wt u_{i+j}v++n+j+m_2.$ It follows from
Lemma \ref{l2.2} that $\Res_{z_2}Y(u_{i+j}v,z_2)w\frac{(1+z_2)^{\wt u+\wt v+n+m_2-1-i}}{z_2^{2n+m_1+m_2+2+j}}$
lies in $\O_n(M),$ as desired.

(2) By Lemma \ref{l2.2}, $\O_n(M)\subset \O_{n-1}(M).$ So it is enough to show
that $u*_nw\equiv u*_{n-1}w, w*_nu\equiv w*_{n-1}u$ modulo  $\O_{n-1}(M)$ for $u\in V$ and $w\in M.$
The proof  is similar to that of Proposition 2.4 of \cite{DLM3}.

(3) We fist prove that $\phi(\O_n(M))\subset \O_n(M).$  Recall the following conjugation formulas from \cite{FHL}:
$$z^{L(0)}Y(u,z_{0})z^{-L(0)}=Y(z^{L(0)}u,zz_{0}),$$
$$e^{zL(1)}Y(u,z_{0})e^{-zL(1)}=Y\left(e^{z(1-zz_{0})L(1)}(1-zz_{0})^{-2L(0)}u,{z_0\over 1-zz_{0}}\right)$$
on $M$ for $u\in V.$ Then for $u\in V$ and $w\in M,$
\begin{eqnarray*}
& &\phi(u\circ_n w)=e^{L(1)}e^{\pi i L(0)}\Res_{z}\frac{(1+z)^{{\wt}u+n}}{z^{2n+2}}Y(u,z)w
\nonumber\\
& &=\Res_{z}\frac{(1+z)^{{\wt}u+n}}{z^{2n+2}}e^{L(1)}
Y((-1)^{L(0)}u,-z)e^{\pi iL(0)}w\nonumber\\
& &=\Res_{z}\frac{(1+z)^{{\wt}u+n}}{z^{2n+2}}Y\left(e^{(1+z)L(1)}(1+z)^{-2L(0)}(-1)^{L(0)}u,{-z\over 1+z}\right)e^{L(1)}e^{\pi iL(0)}w.
\end{eqnarray*}
Making change of variable $z=\displaystyle{-{z_{0}\over 1+z_{0}}}$ and
using the residue formula for the change of variable \cite{Z}
$$\Res_zg(z)=\Res_{z_0}(g(f(z_0)){d\over{dz_0}}f(z_0))$$
(for $g(z)=\sum_{n\geq N}v_nz^{n+\alpha}$ and $f(z)=\sum_{n>0}c_nz^n$
with $n\in\Z$, $v_n\in V,$ $\alpha, c_n\in\C$)
yields
\begin{eqnarray*}
& &\ \ \ \ \phi(u\circ_nw)\\
& &=-\Res_{z_{0}}\frac{(1+z_0)^{-{\wt}u+n}}{z_0^{2n+2}}Y\left(e^{(1+z_{0})^{-1}L(1)}(1+z_{0})^{2L(0)}(-1)^{L(0)}u,z_{0}\right)
 e^{L(1)}e^{\pi iL(0)}w\\
&&=(-1)^{{\wt}u+1}\Res_{z}\frac{(1+z)^{{\wt}u+n}}{z^{2n+2}}Y(e^{(1+z)^{-1}L(1)}u,z)e^{L(1)}e^{\pi iL(0)}w\nonumber\\
& &=(-1)^{{\wt}u+1}\sum_{j=0}^{\infty}{1\over j!}\Res_{z}
\frac{(1+z)^{{\wt}u+n-j}}{z^{2n+2}}Y(L(1)^{j}u,z)
e^{L(1)}e^{\pi iL(0)}w\\
& &=(-1)^{{\wt}u+1}\sum_{j=0}^{\infty}{1\over j!}\Res_{z}
\frac{(1+z)^{{\wt}(L(1)^ju)+n}}{z^{2n+2}}Y(L(1)^{j}u,z)
e^{L(1)}e^{\pi iL(0)}w
\end{eqnarray*}
which lies in $\O_n(M)$ by definition.

Since the proof $\phi(u*_nw)=\phi(w)*_n\phi(u)$ and $\phi(w*_nu)=\phi(u)*_n\phi(w)$ are similar, we give a proof of $\phi(w*_nu)=\phi(u)*_n\phi(w)$ only. We have
\begin{eqnarray*}
&&\phi(w*_nu)=\phi\left(\sum_{m=0}^{n}(-1)^n{m+n\choose n}\Res_zY(u,z)w\frac{(1+z)^{\wt\,u+m-1}}{z^{n+m+1}}\right)\\
& &=\sum_{m=0}^{n}(-1)^n{m+n\choose n}\Res_z\frac{(1+z)^{\wt\,u+m-1}}{z^{n+m+1}}e^{L(1)}Y((-1)^{L(0)}u,-z)e^{\pi iL(0)}w\\
& &=\sum_{m=0}^{n}(-1)^n{m+n\choose n}\Res_z\frac{(1+z)^{\wt\,u+m-1}}{z^{n+m+1}}\\
& &\ \ \ \ \cdot
Y(e^{(1+z)L(1)}(1+z)^{-2L(0)}(-1)^{L(0)}u,\frac{-z}{1+z})e^{L(1)}e^{\pi iL(0)}w
\\
& &=\sum_{m=0}^{n}(-1)^{\wt u+m}{m+n\choose n}\Res_z\frac{(1+z)^{\wt\,u+n}}{z^{n+m+1}}Y(e^{(1+z)^{-1}L(1)}u,z)e^{L(1)}e^{\pi iL(0)}w\\
& &=\sum_{j=0}^{\infty}\frac{1}{
j!}\sum_{m=0}^{n}(-1)^{\wt u+m}{m+n\choose n}\Res_z\frac{(1+z)^{\wt\,u+n-j}}{z^{n+m+1}}
Y(L(1)^{j}u,z)e^{L(1)}e^{\pi iL(0)}w\\
& &=\sum_{j=0}^{\infty}\frac{1}{
j!}(L(1)^j(-1)^{L(0)}u)*_n\phi(w)\\
& &=\phi(u)*_n\phi(w).
\end{eqnarray*}

The proof of (4) is similar to that of Theorem \ref{tha} (3).
\qed

We now take $M=V.$ Theorems \ref{tha} and \ref{t3.5} give the following corollary.
\begin{cor} Let $n\geq 0.$ Then

(1) $\AA_n(V)$ is an associative algebra with the product $*_n$ and the identity $\1+\O_n(V).$

(2) $O_n(V)/\O_n(V)$ is a two sided ideal of $\AA_n(V).$

(3) If $V$ is rational then both $\O_{s-1}(V)/\O_{s}(V)$ and $\AA_{s-1}(V)$ are two sided ideals of $\AA_n(V)$ for $s=1,...,n$ and
$$\AA_n(V)=\AA_0(V)\bigoplus \bigoplus_{s=1}^n\O_{s-1}(V)/\O_{s}(V).$$
\end{cor}

We remark that $\omega+\O_n(V)$ does not lie in the center of $\AA_n(V)$ as $\omega*_nu-u*_n\omega
=(L(-1)+L(0))u$ for $u\in V.$

\section{$\AA_n(M)$ and intertwining operators}

In this section we discuss how $\AA_n(M)$ is related to the intertwining operators and fusion rules. 

We first review the intertwining operators and tensor product of modules from \cite{FHL}.
Let $(W^{i},Y)$ be weak $V$-modules for
$i=1,2,3$. An intertwining operator of type \hspace{-0.2
cm}\singlespace $\left(\!\hspace{-3 pt}\begin{array}{c} W^3\\
W^1 W^2\end{array}\hspace{-3 pt}\!\right)$\doublespace is a linear
map
\begin{eqnarray*}
W^{1}& &\rightarrow ({\rm Hom}(W^{2},W^{3}))\{z\}\\
u & &\mapsto I(u,z)=\sum_{n\in\mathbb{C}}u_{n}z^{-n-1}
\end{eqnarray*}
satisfying the following conditions:

(i)  For any fixed $n\in \mathbb{C}, u\in W^{1}, v\in W^{2}, u_{n+k}v=0$ for sufficiently
large integer $k;$

(ii) $I(L(-1))u,z)=\frac{d}{dz}I(u,z)$ for $u\in
W^{1};$

(iii) The Jacobi identity holds:
\begin{eqnarray*}
&& z_{0}^{-1}\delta
(\frac{z_{1}-z_{2}}{z_{0}})Y(u,z_{1})I(v,z_{2})w-
z_{0}^{-1}\delta (\frac{z_{2}-z_{1}}{-z_{0}})I(v,z_{2})Y(u,z_{1})w\\
&=& z_{2}^{-1}\delta
(\frac{z_{1}-z_{0}}{z_{2}})I(Y(u,z_{0})v,z_{2})w
\end{eqnarray*}
for any $u\in V$, $v\in W^{1}$ and $w\in W^{2}$.

Denote by ${\cal I}\left (\hspace{-3 pt}\begin{array}{c} W^3\\
W^1\,W^2\end{array}\hspace{-3 pt}\right)$ the vector space of all intertwining operators of this
type. We call
$$N_{W^{1}\,W^{2}}^{W^{3}}=\dim\ \hspace{-0.2
cm}\singlespace {\cal I}\left (\hspace{-3 pt}\begin{array}{c} W^3\\
W^1\,W^2\end{array}\hspace{-3 pt}\right)$$
the fusion rules.

Now let $W^{i}=\bigoplus_{n=0}^{\infty}W^{i}(n)$ $(i=1,2,3)$
be three admissible $V$-modules such that $L(0)|_{W^{i}(n)}=(n+h_{i}){\rm id}$
for some constant $h_{i}\in \mathbb{C}$ for  $i=1,2,3$ with $W^{i}(0)\ne 0.$
We write $\deg w=n$ for $w\in W^{i}(n)$. Then
an intertwining operator $I \in {\cal I}\left (\hspace{-3 pt}\begin{array}{c} W^3\\
W^1\,W^2\end{array}\hspace{-3 pt}\right)$
can be written as
\begin{eqnarray*}
I(v,z)=\sum_{n\in\mathbb{Z}}v(n)z^{-n-1}z^{-h_{1}-h_{2}+h_{3}}
\end{eqnarray*}
for $v\in W^1$ (cf. Proposition 1.5.1 in \cite{FZ}). It is clear from the definition that
for homogeneous  $u\in W^{1}, m,n\in \mathbb{Z}$,
\begin{eqnarray*}
u(n)W^{2}(m)\subseteq W^{3}(\deg u+m-n-1).
\end{eqnarray*}

We now can define the tensor product of two admissible $V$-modules. The  tensor product  for the admissible modules $M,N$ is an admissible module $M\boxtimes N$ together with an intertwining operator
$F\in {\cal I}\left(\hspace{-3 pt}\begin{array}{c} M\boxtimes N\\
M\,N\end{array}\hspace{-3 pt}\right)$\doublespace
satisfying the following universal mapping property:
For any admissible $V$-module $W$ and any
intertwining operator $I\in {\cal I}\left(\hspace{-3 pt}\begin{array}{c} W\\
M\,N\end{array}\hspace{-3 pt}\right),$\doublespace there exists a unique
$V$-homomorphism $\psi$ from $M\boxtimes N$ to $W$ such that $I=\psi\circ F.$
Note that the definition of the tensor product does not guarantee the existence of tensor product.

From now on we assume that $V$ is rational and $C_2$-cofinite. Let $M^0,...,M^p$ be the irreducible $V$-modules as before. For short we set $N_{ij}^k=N_{M^i\!M^j}^{M^k}$ for $i,j,k\in\{0,...,p\}.$ Then $N_{ij}^k$
are finite (cf. \cite{ABD}). We also set ${\cal I}_{ij}^k={\cal I}\left (\hspace{-3 pt}\begin{array}{c} W^k\\
W^i\,W^j\end{array}\hspace{-3 pt}\right).$ Let $I_{ij}^{k,s}$ for $s=1,...,N_{ij}^k$  be a basis
of ${\cal I}_{ij}^k.$  The following result is well known:
\begin{thm}\label{tensor} If $V$ is rational then the tensor product $M^i\boxtimes M^j$ exists for any $i,j.$
In fact,
$$M^{i}\boxtimes M^j=\bigoplus_{k=0}^p\bigoplus_{s=1}^{N_{ij}^k}M^{k,s}$$
where $M^{k,s}\cong M^k$ as $V$-module and $F=\sum_{k=0}^p\sum_{s=1}^{N_{ij}^k}I_{ij}^{k,s}$
such that  $I_{ij}^{k,s}(u,z)M^j\subset M^{k,s}\{z\}.$
\end{thm}

The associativity of the tensor product is established in \cite{H} with some additional assumptions:
\begin{thm}\label{huang} If $V$ is rational, $C_2$-cofinite and self-dual, then the tensor product
of $V$-modules is associative.
\end{thm}

We now turn our attention to the connection between intertwining operators and $\AA_n(M).$ Let
$I\in {\cal I}_{ij}^k.$ Then for $w\in M^i$
\begin{eqnarray}\label{e1}
I(w,z)=\sum_{m\in\mathbb{Z}}w(m)z^{-m-1}z^{-\l_{i}-\l_{j}+\l_{k}}
\end{eqnarray}
and $w(\deg w-1-t+s)M^j(s)\subset M^k(t)$ for all $s,t.$ For short we set
$$o_{t,s}^I(w)=w(\deg w-1-t+s)$$
 for homogeneous $w$ and extend it linearly to entire $M^i.$ The following result is a generalization of Lemma 1.5.2 of \cite{FZ} and Theorem 3.2 of \cite{DLM3}.
\begin{lem}\label{ll} The map $(w^i,w^j)\mapsto o^I_{t,s}(w^i)w^j$ for $w^i\in M^i, w^j\in M^j(s)$ for $s,t\leq n$ induces
an $A_n(V)$-module homomorphism $I_{t,s}$ from $\AA_n(M^i)\otimes_{A_n(V)}M^j(s)$ to $M^k(t).$
\end{lem}

\pf As in \cite{FZ} we need to verify that $o^I_{t,s}(u*_nw^i)=o(u)o^I_{t,s}(w^i),$ $o^I_{t,s}(w^i*_nu)=o^I_{t,s}(w^i)o(u)$ and
$o^I_{t,s}(w)=0$ on $M^j(s)$ with $s\leq n$ for $w^i\in M^i$ and $w\in \O_n(M^i).$ The proof is similar to that of Lemma 4.1 of \cite{DJ}.
\qed

The proof of the following theorem is similar to that of Theorem 1.5.2 of \cite{FZ} (see the proof of Theorem 2.11 of \cite{L2}).
\begin{thm}\label{tn1} Assume that $M^j(s),M^k(t)\ne 0$ where $s,t\leq n.$  Then the map from ${\cal I}_{ij}^k$ to
$$\Hom_{A_n(V)}(\AA_n(M^i)\otimes_{A_n(V)}M^j(s), M^k(t))$$
 by sending $I$ to $I_{t,s}$
is a linear isomorphism.
\end{thm}

We need to review a well known result about finite dimensional simple algebra. Let $A$ be a finite dimensional semisimple associative algebra and $W$
an $A$-module. Then $W^*$ is a right $A$-module such that $(fa)(w)=f(aw)$ for $a\in A,$ $f\in W^*$ and
$w\in W.$
\begin{lem}\label{la} If $S,T$ are two simple $A$-modules then $S^*\otimes_{A}T=0$ if $S$ and $T$ are inequivalent, and $S^*\otimes_AS=\C.$
\end{lem}

The following corollary tells us the bimodule structure of $\AA_n(M^i)$ explicitly.
\begin{cor}\label{cn1} The $A_n(V)$-bimodule $\AA_n(M^i)$ has the decomposition
$$\bigoplus_{j,k=0}^p\bigoplus_{s,t=0}^nN_{ij}^kM^k(t)\otimes_{\C}M^{j^*}(s).$$
In particular,
$$\AA_n(V)=\bigoplus_{j=0}^p\bigoplus_{s,t=0}^nM^j(t)\otimes_{\C}M^{j^*}(s).$$
\end{cor}

\pf By Theorems \ref{tha}, \ref{t3.5} and Lemma \ref{la} we know that
$$\AA_n(M^i)=\bigoplus_{j,k=0}^p\bigoplus_{s,t=0}^na_{ij}^{k}(s,t)M^k(t)\otimes_{\C}M^{j^*}(s)$$
as an $A_n(V)$-bimodule for some nonnegative integers $a_{ij}^{k}(s,t).$ Theorem \ref{tn1} then asserts that
$a_{ij}^{k}(s,t)=N_{ij}^k.$ 
\qed

We remark that from the associativity of the tensor product (Theorem \ref{huang}) one can easily prove that
$A(M\boxtimes N)=A(M)\otimes_{A(V)}A(N)$ for any $V$-modules $M,N.$ That is, $A$ is functor from the $V$-module category to the $A(V)$-bimodule category  preserving the tensor product. However, the functor $\AA_n$ does not preserve the tensor product by Corollary \ref{cn1}. 

The following corollary is a refinement of Theorem \ref{t3.5} (4).
\begin{cor} Let $n\geq 1.$ Then as a $A_n(V)$-bimodule
$$\O_{n-1}(M^i)/\O_n(M^i)=\bigoplus_{j,k=0}^p\bigoplus_{s=0}^nN_{ij}^k(M^k(n)\otimes M^{j^*}(s)\oplus M^k(s)\otimes M^{j^*}(n)).$$
\end{cor}

We next investigate the relation among $\AA_n(V)$ and $A_{t,s}(V)$ defined in \cite{DJ}. The $A_{t,s}(V)$
is an $A_t(V)-A_s(V)$-bimodule and the construction of $A_{t,s}(V)$ is motivated from the representation theory of vertex operator algebra. The $A_{t,s}(V)$ is defined as $V/O_{t,s}(V)$ where $O_{t,s}(V)$ is a subspace of $V$ mainly containing  $L(-1)u+(L(0)+s-t)u$ and
$$u\circ_{s}^{t}v={\rm
Res}_{z}\frac{(1+z)^{wtu+s}}{z^{t+s+2}}Y(u,z)v$$
for $u, v\in V.$ One can easily show that
$${\rm
Res}_{z}\frac{(1+z)^{wtu+s+a}}{z^{t+s+2+b}}Y(u,z)v\in O_{t,s}(V)$$
for any nonnegative integer $a\leq b$ (cf. \cite{DJ}).  It is clear from the definition that $\O_n(V)$ is a subspace of $O_{t,s}(V)$ for any $s,t\leq n.$ Using the proof of Theorem \ref{t3.5} we see that
the identity map on $V$ induces an epimorphism of $A_n(V)$-bimodules from $\AA_n(V)$ to $A_{t,s}(V)$ for $s,t\leq n.$ Since $V$ is rational here, $A_{t,s}(V)$ is a sub $A_n(V)$-bimodule of  $\AA_n(V).$ In this case,
$$A_{t,s}(V)=\bigoplus_{i=0}^{p} \bigoplus_{l=0}^{{\rm
min}\{t,s\}}M^i(t-l)\otimes M^{i^*}(s-l)
$$
by Theorem 4.16 of \cite{DJ}.

Finally, we  interprete Corollary \ref{cn1}  in terms of intertwining operators.
Recall a well known fact about associative algebra. If $A$ is an associative algebra and $M,N$ are $A$-modules then $\Hom_{\C}(M,N)$ is an $A$-bimodule defined as follows: for $f\in  \Hom_{\C}(M,N),$ $w\in M$ and $a,b\in A,$ $(afb)(w)=af(bw).$
If both $M,N$ are finite dimensional then $\Hom_{\C}(M,N)$ is isomorphic to $N\otimes_{\C} M^*$ as
$A$-bimodule where $M^*$ has a right $A$-module structure as mentioned before.

We now fix an intertwining operator $I\in{\cal I}_{ij}^k.$ Recall the expression of $I(w,z)$ from equation (\ref{e1}).
Set
$$O^I_{t,s}=\{o^I_{t,s}(w)|w\in M^i\}\subset \Hom_{\C}(M^j(s),M^k(t)).$$
\begin{lem}\label{ll1} Let $s,t\leq n.$ If $I\ne 0$ then $O^I_{t,s}=\Hom_{\C}(M^j(s),M^k(t)).$
\end{lem}
\pf  Note that $o(a)o^I_{t,s}(w)o(b)=o^I_{t,s}(a*_nw*_nb)$ on $M^j(s)$ for $a,b\in V$ and $w\in M^i$
from the proof of Lemma \ref{ll}. It follows immediately that $O^I_{t,s}$ is an $A_n(V)$ sub-bimodule of $\Hom_{\C}(M^i(s),M^k(t)).$ So it is sufficient to prove that $O^I_{t,s}\ne 0$ as long as $\Hom_{\C}(M^i(s),M^k(t))\ne 0.$  Let $0\ne u\in M^j(s).$ Using a  result from
 \cite{DM1} and \cite{L1} we know that
$$M^k=\<w(m)u|w\in M^i,m\in\Z\>.$$
This implies that $M^k(t)$ is spanned by $o^I_{t,s}(w)u$ for $w\in M^i.$ In particular, $O^I_{t,s}\ne 0,$ as expected.
\qed

Recall that $I_{ij}^{k,m}$ for $m=1,...,N_{ij}^k$  is a basis
of ${\cal I}_{ij}^k$ and  $M^{k,m}\cong M^k$ are $V$-modules for $m=1,...,N^k_{ij}$ as before.
Then by Theorem \ref{tensor} and Corollary \ref{cn1} we have
\begin{cor} As $A_n(V)$-bimodule
$$\AA_n(M^i)=\bigoplus_{j,k=0}^p\bigoplus_{m=1}^{N_{ij}^k}\bigoplus_{s,t=0}^nO^{I_{ij}^{k,m}}_{t,s}.$$
\end{cor}

\section{The dual of $\AA_n(M)$}

In this section we study the dual of $\AA_n(M)$ for a $V$-module $M.$ First note that $\AA_n(M)^*$ is also
an $A_n(V)$-bimodule such that for $a,b\in V,$ $f\in \AA_n(M)^*$ and $w\in M,$ $(afb)(w)=f(b*_nw*_na).$

\begin{lem}\label{dual} Let $A$ be a finite dimensional semisimple associative algebra and $M,N$ be the simple $A$-modules. Then the dual of $A$-module $M\otimes N^*$ is isomorphic to $N\otimes M^*.$
\end{lem}

\pf Let $\{\!u_1,...,u_m\!\}$ and $\{\!v_1,...,v_n\!\}$ be bases of $M$ and $N,$ respectively. Let
$\{u_1^*,...,u_m^*\}$ and $\{v_1^*,...,v_n^*\}$ be the dual bases of $M^*$ and $N^*.$ Let $\{e_{ij}|i,j=1,...,m\}$
be a basis of $\End M\subset A$ such that $e_{ij}u_k=\delta_{j,k}u_i$ for all
$1\leq i,j,k\leq m.$ Similarly, let $\{f_{ij}|i,j=1,...,n\}$ be a basis of $\End N\subset A$ such that
$f_{ij}v_k=\delta_{j,k}v_i$ for $1\leq i,j,k\leq n.$ Note that $\{E_{ij}=u_i\otimes v_j^*|i=1,...,m, j=1,...,n\}$
is a basis of $M\otimes N^*$ and $\{F_{st}=v_s\otimes u_t^*|s=1,...,n, t=1,...,m\}$ is a basis of $N\otimes M^*.$
Let $\{E_{ij}^*|i=1,...,m,j=1,...,n\}$ be the dual basis of  $(M\otimes N^*)^*.$

We claim that $f_{st}E_{ij}^*=\delta_{t,j}E_{is}^*$ and $E_{ij}^*e_{ab}=\delta_{ia}E_{bj}^*.$
Since the proofs of two relations are similar we only give the detail for the first relation. From the definition
we know that
$$(f_{st}E_{ij}^*)(E_{ab})=E_{ij}^*(E_{ab}f_{st})=\delta_{b,s}E_{ij}^*(E_{at})
=\delta_{b,s}\delta_{a,i}\delta_{t,j}$$
for any $a,b.$ That is, $f_{st}E_{ij}^*=\delta_{t,j}E_{is}^*.$

It is clear that the linear map from $(M\otimes N^*)^*$ to $N\otimes M^*$ by sending $E_{ij}^*$ to $F_{ji}$
is an $A$-bimodule isomorphism. \qed

From Lemma \ref{dual} we see that there is a natural paring
$$(\cdot,\cdot): M\otimes N^*\times N\otimes M^*\to \C$$
such that $(u\otimes f,v\otimes g)=f(v)g(u)$ for $u\in M,v\in N, f\in N^*,g\in M^*.$ Moreover the following relation
holds
$$(axb, y)=(x,bya)$$
for $x\in M\otimes N^*, y\in N\otimes M^*, a,b\in A.$

\begin{prop}\label{pdual} Let $V$ be a rational, $C_2$-cofinite vertex operator algebra.
Assume that $M^0,...,M^p$ are the irreducible $V$-modules such that $\lambda_i>0$ if $i\ne 0.$ Then
for any $n\geq 0$ we have $\AA_n(M^i)^*$ and $\AA_n(M^{i^*})$ are isomorphic $A_n(V)$-bimodules.
In particular, $A(M^i)^*$ and $A(M^{i^*})$ are isomorphic $A(V)$-bimodules.
\end{prop}

\pf Recall that the fusion matrix $N(i)=(N_{ij}^k)_{j,k=0}^p.$ By Lemma 5.5 of \cite{DJX} we know that
$N(i^*)$ and the transpose $N(i)^T$ of $N(i)$ are the same. That is, $N_{i^*j}^k=N_{ik}^j$ for all $j,k.$
By Corollary \ref{cn1}, we have
$$\AA_n(M^i)=\bigoplus_{j,k=0}^p\bigoplus_{s,t=0}^nN_{ij}^kM^k(t)\otimes_{\C}M^{j^*}(s),$$
and
$$\AA_n(M^{i^*})=\bigoplus_{j,k=0}^p\bigoplus_{s,t=0}^nN_{i^*k}^jM^j(s)\otimes_{\C}M^{k^*}(t).$$
Lemma \ref{dual} then tells us that $\AA_n(M^i)^*$ and $\AA_n(M^{i^*})$ are isomorphic $A_n(V)$-bimodules.
\qed

\end{document}